\documentclass[12 pt,twoside]{amsart}

\usepackage{amsmath,amsfonts,amsthm,amssymb}

\theoremstyle{definition}

\theoremstyle{remark}

\numberwithin{equation}{section}

\newtheorem*{1.1}{1.1}
\newtheorem*{1.2}{1.2}
\newtheorem*{1.3}{1.3}
\newtheorem*{1.4}{1.4 Remark}

\begin{document}

\title [A group of isometries with non-closed orbits]{A group of isometries with non-closed orbits}

\author[H. Abels]{H. Abels}
\address{Fakult\"{a}t f\"{u}r Mathematik, Universit\"{a}t Bielefeld, Postfach 100131, D-33501 Bielefeld, Germany}
\email{abels@math.uni-bielefeld.de}

\author[A. Manoussos]{A. Manoussos}
\address{Fakult\"{a}t f\"{u}r Mathematik, SFB 701, Universit\"{a}t Bielefeld, Postfach 100131, D-33501 Bielefeld, Germany}
\email{amanouss@math.uni-bielefeld.de}
\thanks{During this research the second author was fully supported by SFB 701 ``Spektrale Strukturen und
Topologische Methoden in der Mathematik" at the University of Bielefeld, Germany. He is grateful for its generosity and hospitality.}

\date{}

\subjclass[2000]{Primary 37B05, 54H15; Secondary 54H20.}

\keywords{Proper action, group of isometries, smooth orbit equivalence relation.}

\begin{abstract}
In this note we give an example of a one-dimensional manifold with two connected components and a complete metric whose group of isometries has an orbit which is not
closed. This answers a question of S. Gao and A. S. Kechris.
\end{abstract}

\maketitle

\section{Preliminaries and the construction of the example}
In \cite[p. 35]{kechris} S. Gao and A. S. Kechris asked the following question. Let $(X,d)$ be a locally compact complete metric space with finitely many
pseudo-components or connected components. Does its group of isometries have closed orbits? This is the case if $X$ is connected since then the group of isometries
acts properly by an old result of van Dantzig  and van der Waerden \cite{d-w} and hence all of its orbits are closed. The above question arose in the following
context. Suppose a locally compact group with a countable base acts on a locally compact space with a countable base. Then the action has locally closed orbits (i.e.
orbits which are open in their closures) if and only if there exists a Borel section for the action (see \cite{glimm}, \cite{effros1}) or, in other terminology, the
corresponding orbit equivalence relation is smooth. For isometric actions it is easy to see that an orbit is locally closed if and only if it is closed. In this note
we give a negative answer to the question of Gao and Kechris. Our space is a one-dimensional manifold with two connected components, one compact isometric to $S^1$,
and one non-compact, the real line with a locally Euclidean metric. It has a complete metric whose group of isometries has non-closed dense orbits on the compact
component. In the course of the construction we give an example of a 2-dimensional manifold with two connected components one compact and one non-compact and a
complete metric whose group $G$ of isometries also has non-closed dense orbits on the compact component. The difference is that $G$ contains a subgroup of index 2
which is isomorphic to $\mathbb{R}$.

Let $(Y,d_1)$ be a metric space. Later on $Y$ will be a torus with a flat Riemannian metric. Let $Z=Y\cup (Y\times\mathbb{R})$. We fix two positive real numbers $R$
and $M$. We endow $Z$ with the following metric $d$ depending on $R$ and $M$.
\[
\begin{split}
&d(y_1,y_2) = d_1(y_1,y_2)\\
&d((y_1,t_1),(y_2,t_2))= d_1(y_1,y_2) + \min (|t_1-t_2|,M)\\
&d(y_1,(y_2,t_2)) = d((y_2,t_2),y_1) = d(y_1,y_2) +R,
\end{split}
\]
for $y_1,y_2\in Y$ and $t_1,t_2\in\mathbb{R}$. It is easy to check that $d$ is a metric on $Z$ if $2R\geq M$. The metric space $Z$ has the following properties

\begin{1.1}
a) For a given point $(y,r)\in Y\times\mathbb{R}$ there is a unique point in $Y$ which is closest to $(y,r)$, namely $y$.

b) Given a point $y\in Y$ the set of points in $Y\times\mathbb{R}$ which are closest to $y$ is the line $\{ y\}\times\mathbb{R}$.

c) For every point $(y,r)\in Y\times\mathbb{R}$ and every $y'\in Y$ there is a unique point on the line $\{ y'\}\times\mathbb{R}$ which is closest to $(y,r)$, namely
$(y',r)$.

d) Let $g_Y$ be an isometry of $Y$ and let $g_{\mathbb{R}}$ be an isometry of the Euclidean line $\mathbb{R}$. Define a map $g=g(g_Y,g_{\mathbb{R}}): Z\to Z$ by
$g|Y:=g_Y$ and $g(y,r)=(g_Y(y),g_{\mathbb{R}}(r))$ for $(y,r)\in Y\times\mathbb{R}$. Then $g$ is an isometry of $Z$.

e) Every isometry of $Z$ is of the form given in d) if $Y$ is compact.
\end{1.1}
\begin{proof}
a) through d) are easily checked. To prove e) let $g$ be an isometry of $Z$. Then $g(Y)=Y$ and $g(Y\times\mathbb{R})=Y\times\mathbb{R}$, since $Y$ is compact and
$Y\times\mathbb{R}$ consists of non-compact components. Then $g_Y:=g|Y$ is an isometry of $Y$. The map $g(g_Y,\mathit{id})^{-1} \circ g$, where $\mathit{id}$ denotes
the identity map, is an isometry of $Z$ which fixes $Y$, hence maps every line $\{ y\}\times\mathbb{R}$ to itself, by b). Let $h_y:\mathbb{R}\to\mathbb{R}$ be defined
by $g(y,t)=(y,h_y(t))$. Then $h_y$ is an isometry of the Euclidean line $\mathbb{R}$ for every $y\in Y$ and all the $h_y$'s are the same, by c), say
$h_y=g_{\mathbb{R}}$. Thus $g=(g_Y,g_{\mathbb{R}})$.
\end{proof}

\begin{1.2}
Let now $Y$ be a 2-dimensional torus with a flat Riemannian metric. $Y$ is also an abelian Lie group whose composition we write as multiplication. Every translation
$L_x$ of $Y$, $L_x(y)= x \cdot y$, is an isometry. Let $g(t)$, $t\in\mathbb{R}$, be a dense one parameter subgroup of $Y$. Let $H\subset Y\times\mathbb{R}$ be its
graph, $H=\{ (g(t),t)\,;\,t\in\mathbb{R} \}$. Our example is $X=Y\cup H$ with the metric induced from $Z=Y\cup (Y\times\mathbb{R})$.
\end{1.2}

\begin{1.3}
a) If $g_{\mathbb{R}}$ is an isometry of the Euclidean line $\mathbb{R}$ then there is a unique isometry $g$ of $X$ such that $g(y,t)\in Y\times \{
g_{\mathbb{R}}(t)\}$. If $g_{\mathbb{R}}$ is the translation by $a$, so $g_{\mathbb{R}}=L_a$ with $L_a(t)=t+a$, then $g$ is the restriction of $g(L_{g(a)},L_a)$ to
$X$. If $g_{\mathbb{R}}$ is the reflection at $O$, $g_{\mathbb{R}}=\mathbf{-1}$, then $g$ is the restriction of $g(\mathit{inv},\mathbf{-1})$ to $X$, where
$\mathit{inv}:Y\to Y$, $\mathit{inv}(y)=y^{-1}$. The reflection in $a\in\mathbb{R}$ is the composition $L_{-2a}\circ (\mathbf{-1})=\mathbf{-1}\circ L_{2a}$.

b) Every isometry of $X$ is of the form in a). It follows that the group of isometries of $X$ has dense non-closed orbits on $Y$ and the other component $H$ is one
orbit.

c) $H$ is locally isometric to the real line with the Euclidean metric, actually $d((g(t),t),(g(s),s))=(1+\| \overset\bullet{g} (0)\|)\, |t-s|$ for small $|t-s|$,
where $\overset\bullet{g} (0)$ is the tangent of the one-parameter group $g(t)$, $t\in\mathbb{R}$, and $\| \cdot\|$ is the norm on the tangent space of $Y$ at the
identity element derived from the Riemannian tensor.
\end{1.3}
\begin{proof}
c) follows from the definition of the metric $d$ on $Y\times\mathbb{R}$. The maps given in a) are isometries of $Z$ and map $X$ to $X$, hence are isometries of $X$.
To prove the uniqueness claim in a) it suffices to prove it for $g_{\mathbb{R}}=\mathit{id}$. But then $g$ is the identity on the image of the one-parameter group
$g(t)$, $t\in\mathbb{R}$, by 1.1 a) and hence on all of $Y$. Hence $g$ has the form given by 1.1 d). To show b) it suffices to show that every isometry $h$ of $H$ is
of the form given in a). This follows from c).
\end{proof}

\begin{1.4}
In our example the space has dimension 2 and the group of orientation preserving isometries is of index 2 in the group of all isometries and is isomorphic to
$\mathbb{R}$. We can reduce the dimension of our space to 1 to obtain a group of isometries with closed orbits on the non-compact component, which is diffeomorphic
and locally isometric to $\mathbb{R}$, and non-closed dense orbits on the compact component, which isometric to $S^1$. The example is as follows. Take a
one-dimensional subtorus $Y_1$ of $Y$ containing the identity element of $Y$. Define $X_1=Y_1\cup H\subset Y\cup H$. Then the group of isometries of $Y_1$ consists of
those maps $g_a=g(L_{g(a)},L_a)$ restricted to $Y_1$ with $g(a)\in Y_1$, and of the maps $g(\mathit{inv}\circ L_{g(2a)},\mathbf{-1}\circ L_{a})$ restricted to $Y_1$
with $g(2a)\in Y_1$. The proof follows from the proof of 1.3.
\end{1.4}

\end{document}